\documentclass[11pt]{elsarticle}
\usepackage{amssymb}
\usepackage{amsmath}
\usepackage{amsthm}
\usepackage{epsfig}
\usepackage{color}
\usepackage{makeidx} 
\usepackage{fancyhdr}

\usepackage{lineno}

\newtheorem{theorem}{Theorem}
\newtheorem{lemma}{Lemma}
\newtheorem{corollary}{Corollary}

\begin{document}

\begin{frontmatter}

\title{Conditional Law and Occupation Times of Two-Sided Sticky Brownian Motion\footnote{This manuscript version is made available under the CC-BY-NC-ND 4.0 license http://creativecommons.org/licenses/by-nc-nd/4.0/ }$^,$\footnote{https://doi.org/10.1016/j.spl.2020.108856}}
\author{Bugra Can\footnote{Current address: Rutgers University, NJ, USA}  and Mine \c{C}a\u{g}lar\footnote{Corresponding author}}
\address{Department of Mathematics, Ko\c{c} University, Istanbul, Turkey}

\begin{abstract}
Sticky Brownian motion on the real line can be obtained as a weak solution of a system of stochastic differential equations. We find the conditional distribution of the process given the driving Brownian motion, both at an independent exponential time and at a fixed time $t>0$. As a classical problem, we find the distribution of the occupation times of a half-line, and at $0$, which is the sticky point for the process.
\end{abstract}


\end{frontmatter}

\section{Introduction}

An important special case of diffusions with boundary conditions  is sticky Brownian motion, defined by means of an infinitesimal generator on the half-line with a sticky boundary condition at 0 \cite{Feller}. Ito and McKean \cite{itoMc} provided a construction of the process as a time-changed reflected Brownian motion.
Ikeda and Watanabe \cite{IkedaWatanabe}  studied such diffusions with boundary conditions as solutions to stochastic differential equations. A sticky  Brownian motion on the  real line, where 0 is not at the boundary, is  characterized as a solution of a system of stochastic differential equations. Recently, both versions have attracted further research interest motivated by applications \cite{Ankir2, Ankir1,bass,Engelbert,hatem,Shckol,Warren}. In this paper, we derive the conditional law of the real-valued sticky Brownian motion at a snapshot in time, given the driving Brownian motion, as well as the occupation times at each half-line and  $0$ based on  the results of \cite{Engelbert,Warren}.

The one-sided, or reflecting, sticky Brownian motion is determined  as a weak solution of  the stochastic differential equation
\begin{equation}
dX_t=1_{\{X_t>0\}}dB_t+\theta 1_{\{X_t=0\}}dt \label{eq:one-sided}
\end{equation}
where $\theta>0$ is the stickiness parameter, with $X_0\ge 0$ \cite[Eq.(1.1)]{Warren}, \cite[Eq.(3.3) and Thm.5]{Engelbert}. In \cite{Warren}, the conditional distribution of $X_t$  given  ${\cal F}_t$, where $({\cal F}_t)_{t\ge 0}$ is the filtration generated by the driving Brownian motion $B$, has been found as
\begin{equation}  \label{2}
\mathbb{P}\{X_t<x|{\cal F}_t\}=e^{-2\theta(B_t+L_t-x)}
\end{equation}
where $L_t=\underset{s\leq t}{\sup}\{-B_s \wedge 0\}$, $x\in [0,B_t+L_t]$, and initial condition $(X_0, B_0)=(0,0)$. This remarkable result is verified from the point of view of stochastic flows in \cite{hatem}.  On the other hand, the two-sided, or global, version of sticky Brownian motion is described by a system of stochastic differential equations given by
\begin{eqnarray}
dX_t&=& 1_{\{X_t \neq 0\}}dB_t \label{eq:sticky1}\\
	1_{\{X_{t}=0\}}dt&=& \frac{1}{\theta}dl_t^{0}(X)\label{eq: sticky2}
\end{eqnarray}
where $l_t^{0}(X)$ denotes the local time of $X$ at $0$. We  will refer to the two-sided version as \emph{the sticky Brownian motion } in this paper. It has been recently shown in \cite{bass} and \cite{Engelbert}  that there exists no strong solution to equations \eqref{eq:sticky1} and \eqref{eq: sticky2}, but the system admits a weak solution uniquely determining the probability law of $X$.
Moreover,  if $(X_t,B_t)$ solves \eqref{eq:sticky1}-\eqref{eq: sticky2}, then $(|X_t|,\hat{B}_t)$ with $\hat{B}_t:= \int_0^t \mbox{sgn}(X_s)dB_s$ solves \eqref{eq:one-sided} as shown in the proof of  \cite[Thm.5,pg.18-19]{Engelbert}. 
As a result, $|X_t|$ is a one-sided sticky Brownian motion with the same parameter $\theta$. 
Based on this association, the positive and negative parts of $X$ behave similar to the one-sided sticky Brownian and play a crucial role in our results.

Our first result is the characterization of the joint distribution of $(X_t,B_t)$  through their  resolvent density as given in Theorem \ref{theo:decompositionXandB}.   This  is then used in Theorem \ref{thm2} to derive the conditional law of $X_T$ given $B_T$ at an independent exponential  time $T$ with parameter $\lambda$ as
\begin{equation}\label{above}
\mathbb{P}\{X_{T}<x|B_T=b\}= \Big( \frac{\mbox{sgn}(b-x)}{4}+\mbox{sgn}(x) \frac{\gamma}{4\delta}\Big)e^{-\delta|b-x|+\gamma(|b|-{sgn}(x)x)}+\frac{i_x(b)}{2}
\end{equation}
where $\gamma^2=2\lambda$, $\delta^2=2\lambda+2\theta\gamma$, and $i_x(b)=1_{\{b<x\}}$.  We invert the Laplace transforms involved in \eqref{above} to obtain $\mathbb{P}\{X_{t}<x|B_t=b\}$ in Corollary \ref{theo: conditionalProbability}. Our approach is parallel to \cite{Warren} by relying on Knight's theorem. However, we condition on a single random variable $B_t$ rather than the filtration ${\cal F}_t$.

As another result, we derive the distributions of the occupation times of the sticky Brownian motion in Theorem \ref{thm3}. Our exposure is along the lines of Brownian motion, for which the celebrated Arcsine law holds \cite{RevuzYor}. The counterparts for the one-sided case are derived based on the results of \cite{Warren} for finding the distribution of the occupation time at 0. The occupation time at 0 during $[0,1]$,  $A_1^0$, is nonzero as the process is sticky at 0, and satisfies
\[
\mathbb{P}\{A^0_1>t\}=\text{Erfc}(\frac{\theta t}{\sqrt{2(1-t)}})
\]
for $0<t<1$. The distribution of the occupation time of the positive axis, $A_1^+$, is found  as
\[
 \mathbb{P}\{A^+_1 > t \}=  1-\mathbb{E}\, \left[ \text{Erfc}
\left( \frac{ -Z^2+Z\sqrt{Z^2+\theta^2(1-t)}} {\theta \sqrt{2t}}\right) \right]
\]
which approximates  Arcsine law when $\theta \rightarrow \infty$ as expected.

The paper is organized as follows. Section 2 gives the resolvents involved and the conditional law for the sticky  Brownian motion. The occupation distributions are found in Section 3.
	
\section{Joint Resolvent and Conditional Distribution}

In this section, we first obtain some fundamental identities using Skorohod
representations, which will also be useful to derive the occupation time distributions in the next section. Then,  we  recover the resolvent of the sticky Brownian motion $X$ found as a weak solution of \eqref{eq:sticky1}-\eqref{eq: sticky2}, with $X_0=0$. As for the main results, we derive the joint resolvent of $X$ and the driving Brownian motion $B$, and the conditional law of $X_t$ given $B_t$ for each $t>0$.

The occupation times of $X_t$  and their right continuous inverses are defined as follows
	\begin{eqnarray}
		A_{t}^{0}=\int_{0}^{t}1_{\{X_{s}=0\}}ds &&\alpha_{t}^{0}=\inf\{u: A_{u}^{0}>t\} \label{eq: zerotime}\\
		A_{t}^{+}=\int_{0}^{t}1_{\{X_{s}>0\}}ds&&\alpha_{t}^{+}=\inf\{u: A_{u}^{+}>t\}\label{eq: positivetime}\\
		A_{t}^{-}=\int_{0}^{t}1_{\{X_{s}<0\}}ds&&\alpha_{t}^{-}=\inf\{u: A_{u}^{-}>t\}\label{eq:negativetime}
	\end{eqnarray}
Note that $A_{t}^{0}$, $A_{t}^{-}$ and $A_{t}^{+}$ are right continuous and non-decreasing processes justifying the definition of the corresponding right inverses.

Let $\mathbb{P}^{(x,a)}$ and $\mathbb{E}^{(x,a)}$ denote the probability law and the expectation, respectively,  given the initial condition $(X_0,B_0)=(x,a)$, $x,a \in \mathbb{R}$. We simply use  $\mathbb{P}$ and $\mathbb{E}$ when $(X_0,B_0)=(0,0)$.
Let  $X^+$ and $X^-$ denote the positive and negative parts of the real-valued random variable $X$, respectively. The identities shown in the following lemma will be used in the sequel.

\begin{lemma} \label{W}
There exist independent Brownian motions $W_{t}^{+}$ and $W_{t}^{-}$ such that
\[ X_{\alpha_{t}^{+}}^{+}=W^{+}_{t}+\frac{\theta}{2}A_{\alpha_{t}^{+}}^{0}
\; \; , \quad X_{\alpha_{t}^{-}}^{-}=W^{-}_{t}+\frac{\theta}{2}A_{\alpha_{t}^{-}}^{0}
\]
and
\[
A^{0}_{\alpha_{t}^{+}}=\frac{2}{\theta}\, \underset{s\leq t}{\sup}(-W^{+}_{s}) \; \; , \quad A^{0}_{\alpha_{t}^{-}}=\frac{2}{\theta}\, \underset{s\leq t}{\sup}(-W^{-}_{s}) \; .
\]
\end{lemma}

\begin{proof}
Tanaka formula \cite[Thm.VI.1.2]{RevuzYor} for $X$ satisfying equations \eqref{eq:sticky1}-\eqref{eq: sticky2} yields
\begin{eqnarray}
	X_{t}^{+}&=&\int_{0}^{t}1_{\{X_{s}>0\}}dB_{s}+\frac{\theta}{2}A_{t}^{0} \nonumber  \\
	X_{t}^{-}&=& - \int_{0}^{t}1_{\{X_{s}< 0\}}dB_{s}+\frac{\theta}{2}A_{t}^{0} \; .  \label{eq:Tanaka2}
\end{eqnarray}
As a result of a time change $t \rightarrow \alpha_t^+$, we get
\[
	X_{\alpha_{t}^{+}}^{+}= \int_{0}^{\alpha_{t}^{+}}1_{\{X_{s}>0\}}dB_{s}+\frac{\theta}{2}A_{\alpha_{t}^{+}}^{0}\label{eq:distA0}
\]
Define
\[
W^+_t:=\int_{0}^{\alpha_{t}^{+}}1_{\{X_{s}>0\}}dB_{s} \quad \mbox{and} \quad W^-_t:=- \int_{0}^{\alpha_{t}^{-}}1_{\{X_{s}<0\}}dB_{s}.
\]
Then, $\langle W^+_t,W^+_t\rangle=\langle W^-_t,W^-_t\rangle=t$    by definitions \eqref{eq: positivetime}  and \eqref{eq:negativetime}. This implies that $W^+$ and $W^-$ are Brownian motions by Levy's characterization theorem.
Note that   $\langle W^{+}_{A^+_t},W^+_{A^+_t}\rangle=A^+_t$, $\langle W^{-}_{A^-_t},W^+_{A^-_t}\rangle=A^-_t$ and $\langle W^-_{A^-_t},W^+_{A^+_t}\rangle=0$. Therefore, by Knight's theorem (\cite[Thm.V.1.9]{RevuzYor}) $W^+_t $ and $W^-_t $ are independent Brownian motions as $A^+_\infty=A^-_\infty =\infty $.
 Furthermore, $\frac{\theta}{2}A^0_{\alpha_{t}^{+}}$ is a continuous and strictly increasing function of $\alpha_{t}^{+}$ and constant on the set $\{t: X^{+}_{\alpha_{t}^{+}} > 0 \}$.
 Thus, the result follows by Skorohod's Lemma \cite[Lem.VI.2.1]{RevuzYor}.
\end{proof}

The following lemma gives the  resolvent kernel of the sticky Brownian motion $X$ starting from zero. Although it is available in \cite[pg.23]{Howitt}, where $\alpha^0$ and  analogous definition for the right continuous inverse of the time spent in $\mathbb{R} \setminus \{0\}$   are considered, we provide a proof referring to the times spent by $X$ in the three regions introduced above inspired by the  proof of \cite[Prop.13]{Warren} for one-sided sticky. Our approach emphasizes the two-sided nature of $X$, also used in Section 4.

\begin{lemma} \label{lem2} The resolvent kernel $p_{\lambda}(0,dy)$  of the sticky
Brownian motion $X$ starting from 0 is given as
\[
p_{\lambda}(0,dy)=\frac{\theta }{\theta \gamma+\lambda }e^{-\gamma |y|}dy+\frac{1}{\theta\gamma+\lambda}\delta_{0}(dy)
\]
where $\gamma^2=2\lambda$. 
\end{lemma}

\begin{proof}
Take exponentially distributed random times $T_{1}, T_{2} $ and $T_{3}$ with rate $\lambda$ independently from $(X,W)$, and define
\begin{equation}\label{def: exp_T}
T=\alpha_{T_{1}}^{0} \wedge \alpha_{T_{2}}^{+} \wedge \alpha_{T_{3}}^{-}.
\end{equation}
By Lemma \ref{W}, $(\theta/2)A_{\alpha_{t}^{+}}^0$ and $(\theta/2)A_{\alpha_{t}^{-}}^0$ are the running supremums of independent Brownian motions $W^+$ and $W^-$, respectively.
Therefore, $A_{\alpha_{T_2}^{+}}^0$ and $A_{\alpha_{T_3}^{-}}^0$ are exponentially distributed with rate $\frac{\theta \gamma}{2}$, where $\gamma^2=2\lambda$.  Note that $T$ is also exponentially distributed with rate $\lambda$ independent of $X$, similar to the one-sided case given in \cite{Warren}, essentially due to independence of $T_i$, $i=1,\ldots,3$.
Now, we have  $X_{T}=0$ if and only if $T=\alpha_{T_{1}}^{0} $, which can be true only if $\alpha_{T_{1}}^{0}<\alpha_{T_{2}}^{+}\wedge \alpha_{T_{3}}^{-}$. Since $A^{0}_{t}$ is a right continuous and nondecreasing process, if $\alpha_{T_{1}}^{0}<\alpha_{T_{2}}^{+}\wedge \alpha_{T_{3}}^{-}$ then $T_{1}<A^{0}_{\alpha_{T_{2}}^{+}}\wedge A^{0}_{ \alpha_{T_{3}}^{-}}$.  Thus, we get
\begin{eqnarray*}
	\mathbb{P}\{X_{T}=0\}&=&\mathbb{P}\{T=\alpha_{T_{1}}^{0}\}\\
	&=& \mathbb{P}(T_{1}<A^{0}_{\alpha_{T_{2}}^{+}}\wedge A^{0}_{\alpha_{T_{3}}^{-}})\\
	&=&\frac{\lambda}{\lambda+\theta\gamma} \; .
\end{eqnarray*}
For $y\neq 0$, we consider only $y>0$  since on negative axis the calculations are similar. For $X_{T}>0$, we should have $T=\alpha_{T_{2}}^{+}$ and this can be possible only if $A^{0}_{\alpha_{T_{2}}^{+}}<T_{1}\wedge A^{0}_{\alpha_{T_{3}}^{-}}$. Then, note that $X^{+}_{\alpha_{T_{2}}^{+}}$ is a reflected Brownian motion and its probability law for $y\neq 0$ is independent from the level of the running supremum of $-W^+$ (or negative of the running infimum of $W^+$), which is equal to $A^{0}_{\alpha_{T_{2}}^{+}}$ by Lemma 1. This is due to the fact that the reflected Brownian motion is formed by the excursions of $W^+$ from its infimum and the excursion measure does not depend on the level of the last infimum (see e.g. \cite[Sec.XII.2]{RevuzYor}).
Moreover, $X^{+}_{\alpha_{T_{2}}^{+}}$ is independent from $A^{0}_{\alpha_{T_{3}}^{-}}$ again by Lemma \ref{W}. Therefore, $X^+_{\alpha_{T_{2}}^{+}}$ is independent from the event $\{
	A^{0}_{\alpha_{T_{2}}^{+}}<T_{1}\wedge A^{0}_{\alpha_{T_{3}}^{-}}\} $. Then, for $y>0$, we have
\begin{eqnarray}
	\mathbb{P}\{X_{T}\in dy\}&=& \mathbb{P}\{X^+_{\alpha_{T_{2}}^{+}}\in dy,
	A^{0}_{\alpha_{T_{2}}^{+}}<T_{1}\wedge A^{0}_{\alpha_{T_{3}}^{-}}\}\nonumber\\
	&=& \mathbb{P}\{X^+_{\alpha_{T_{2}}^{+}}\in dy\}\mathbb{P}\{
	A^{0}_{\alpha_{T_{2}}^{+}}<T_{1}\wedge A^{0}_{\alpha_{T_{3}}^{-}}\}\nonumber \\
	&=&   \gamma e^{-\gamma y} \; \frac{\frac{\theta \gamma}{2}}{\lambda+ \theta \gamma } \, dy\; = \;\frac{\theta \lambda }{\lambda+ \theta \gamma }e^{-\gamma |y|}dy \label{prob: Xiny}
\end{eqnarray}
where we have used the fact that the distribution of the reflected Brownian motion at $T_2$ is  exponential with parameter $\gamma=\sqrt{2\lambda}$.
For $y<0$, we have the same expression due to symmetry. The result follows as $ p_\lambda(0,dy)= (1/\lambda)\mathbb{P}(X_{T}\in dy)$.
\end{proof}

We are ready to find the joint resolvent $\mathcal{V}_{\lambda}$ of the process $(X,B)$, defined by 
\[
\mathcal{V}_{\lambda}f(x,a) = (1/\lambda) \mathbb{E}^{(x,a)}[f(X_T,B_T)]
\]
for bounded measurable functions $f:\mathbb{R}^2\rightarrow \mathbb{R}$, where $T$ is an independent exponential time with parameter $\lambda$.  

\begin{theorem}\label{theo:decompositionXandB}
Let $(X_t,B_t)$ be the solution to stochastic differential equations \eqref{eq:sticky1} - \eqref{eq: sticky2} with $X_0=0$. Then, there exists a Brownian motion $W^0_t$   such that
\begin{equation}
X_t+W^0_{A_{t}^0}=B_t \label{eq:decompositionXandB}
\end{equation}
and the joint resolvent $\mathcal{V}_{\lambda}$ of the process $(X_t,B_t)$ is given by
 \begin{eqnarray*}
\mathcal{V}_{\lambda}f(x,a)&=&\frac{1}{\theta}\int_{-\infty}^{\infty}\nu_{\lambda}
(x,a,0,b)f(0,b)\,db+\int_{-\infty}^{\infty}\int_{-\infty}^{\infty}\nu_{\lambda}
(x,a,y,b)f(y,b)\,dy\,db\\
&& \quad +\int_{\{a+y=b+x\}}r_{\lambda}^{-}(x,y)f(y,b)dy
\end{eqnarray*}
where
$$
\nu_{\lambda}(x,a,y,b)=\frac{\theta}{2\delta}e^{-\gamma(|y|+|x|)-\delta |b-y-a+x|}
$$
with $\gamma^2=2\lambda$, $\delta^2=2\lambda+2\theta\gamma$, and
$$r_{\lambda}^{-}(x,y)=\gamma^{-1}(e^{-\gamma|y-x|} -e^{-\gamma|y+x|})
$$
which is the resolvent density of Brownian motion killed at 0.
\end{theorem}

 \begin{proof}
When $X_0=0$, \eqref{eq:sticky1} - \eqref{eq: sticky2}  imply
\[
	X_{t}+\int_{0}^{t}1_{\{X_{s}=0\}}dB_{s}=B_{t}.
\]
Define $W^0_{t}$ as
\begin{equation*}
	W^0_{t}=\int_{0}^{\alpha_{t}^{0}}1_{\{X_{s}=0\}}dB_{s}.
\end{equation*}
Then, Levy's characterization and definition of $\alpha_{t}^{0}$ yield $W^0_{t}$ is a Brownian motion and so decomposition \eqref{eq:decompositionXandB} is shown.

The resolvent of $(X,B)$ can be derived by similar methods used for the resolvent kernel of $X$ as in Lemma \ref{lem2}.
Recall the definition of $T$ in  \eqref{def: exp_T} as $T=\alpha_{T_{1}}^{0} \wedge \alpha_{T_{2}}^{+} \wedge \alpha_{T_{3}}^{-}$. Now, we have
 \begin{equation} \label{A0}
  A^0_T=T_1 \wedge A^0_{\alpha^+_{T_2}} \wedge A^0_{\alpha^{-}_{T_3}}
  \end{equation}
where $A^0_{\alpha^+_{T_2}},A^0_{\alpha^-_{T_3}} $ and $T_1$ are independent exponential random variables with rates $\frac{\theta\gamma}{2},\frac{\theta\gamma}{2}$ and $\lambda$, respectively. The distributions of $A^0_{\alpha^+_{T_2}}$ and $A^0_{\alpha^-_{T_3}}$, together with their independence follow from  Lemma \ref{W}. Therefore, $A^0_T$ is  exponentially distributed with  rate $\lambda+\theta\gamma$.   Moreover, $A^0_{\alpha^{+}_{T_2}}$, $A^0_{\alpha^{-}_{T_3}}$, and $W^0$  are  all independent because   $W_t^+$ and  $W_t^-$, and $W_t^0$ are independent Brownian motions again by Knight's theorem. It follows that the $W^0$ is independent of $A^0_T$ given in \eqref{A0}, and $W^0_{A^0_T}$ is double exponentially distributed with parameter $\delta:=\sqrt{2(\lambda+\theta\gamma)}$, with the probability density function $f(x)=(\delta/2)e^{-\delta |x|}$ for $x\in \mathbb{R}$. 

To find the joint resolvent, we first consider a starting point of the form $(X_0,B_0)=(0,a)$. Then,  \eqref{eq:decompositionXandB} becomes
\[
  X_{t}+a+ W^0_{A^0_t}=B_{t}
\]
We will consider the events $X_T=0$ and $ X_T \neq 0 $ separately. The first event  corresponds to the event  $T=\alpha^0_{T_1}$, whereas the second one occurs only if  $T=\alpha^+_{T_2}$  or $\alpha^-_{T_3}$. Since  $T=\alpha^0_{T_1}$ is equivalent to $A^0_T=T_1$, we have
\begin{eqnarray*}
\mathbb{P}^{(0,a)}\{X_{T}=0, B_{T}\in db \}
&=& \mathbb{P}^{(0,a)}\{X_{T}=0, W^0_{A^0_T}\in d(b-a) \}\\
&=& \mathbb{P}^{(0,a)}\{ A^0_T=T_1, W^0_{A^0_T}\in d(b-a) \}\\
&=&\frac{\lambda}{2\delta}e^{-\delta |b-a|}db
\end{eqnarray*}
by independence of $A^0_T$ and $W^0$ as explained after \eqref{A0} above, and $T_1$ being the minimum of three independent exponential random variables.
On the other hand, when $ X_T \neq 0 $, it is implicit that $T=\alpha^+_{T_2}$  or $\alpha^-_{T_3}$.   
Then, for $y>0$ ($y<0$), we have  
\begin{eqnarray*}
\mathbb{P}^{(0,a)}\{X_{T}\in dy, B_{T}\in db \}
	&=&\mathbb{P}^{(0,a)}\{ X_{T}\in dy,W^0_{A_{T}^{0}} \in d(b-y-a) \}\\ 
	&=&\mathbb{P}^{(0,a)}\{X_{T}\in dy\}\, \mathbb{P}^{(0,a)} \{ W^0_{A_{T}^{0}} \in d(b-y-a) \}\\ 
	&=&\frac{\lambda}{2\delta}e^{-\gamma|y|-\delta |b-y-a|}\,dy\,db
\end{eqnarray*}
where we have used \eqref{prob: Xiny}, and the independence of $X_T$  and $W^0_{A_{T}^{0}}$ when $T=\alpha^+_{T_2}$  ($\alpha^-_{T_3}$). To see the independence, recall that $(W^+_{ t}, W^-_{t},W^0_t)$ is a 3-dimensional Brownian motion by Knight's theorem. Therefore, $X^+_{\alpha^+_\cdot}$, $ X^-_{\alpha^-_\cdot}$ and $W^0$ are  independent processes by Lemma \ref{W}. Moreover, $X^+_{\alpha_{T_{2}}^{+}}$ is independent from the event $\{	A^{0}_{\alpha_{T_{2}}^{+}}<T_{1}\wedge A^{0}_{\alpha_{T_{3}}^{-}}\} $ as explained in the proof of Lemma \ref{W} before \eqref{prob: Xiny}, and these together characterize the distribution of $X_T$. It follows that $X_T$ is independent from $W^0_{A_{T}^{0}}$, analogous to the one-sided case given in   \cite[pg.12]{Warren}. 

Now, as $\mathcal{V}_{\lambda}f(x,a) = (1/\lambda) \mathbb{E}^{(x,a)}[f(X_T,B_T)]$, we get
\begin{eqnarray*}
\mathcal{V}_{\lambda}f(0,a)& = & \int \int \mathbb{P}^{(0,a)}\{X_{T}\in dy, B_{T}\in db \}f(y,b) dy\, db  \\
&& \; + \int \mathbb{P}^{(0,a)}\{X_{T}=0, B_{T}\in db \} f(0,b) \, db \\
&=&\frac{1}{\theta} \int \int \nu_{\lambda}
(0,a,y,b)f(y,b)\,dy\,db+  \int \nu_{\lambda}
(0,a,0,b)f(0,b)\,db
\end{eqnarray*}
where $\nu_{\lambda}(x,a,y,b):=\frac{\theta}{2\delta}e^{-\gamma(|y|+|x|)-\delta |b-y-a+x|}$ is used for brevity of notation. 

To find the resolvent of $(X_t,B_t)$ starting from $(x,a)$, $x\neq 0$, let us define the time $H_0$ at which $X_t$ hits zero for the first time, and use the strong Markov property of $X_{t}$ at $H_0$, the hitting time of $X$ to 0.  Since $X_t$ is in lockstep with $B_t$ for   $t\in [0, H_{0}]$, having the same trajectory   according to \eqref{eq:sticky1} - \eqref{eq: sticky2}, the law of the process $X$ during $t\in [0, H_{0}]$ is just that of a Brownian motion started at $x$ and killed at $0$. If we define $\mathcal{R}_{\lambda}^{-}$ as the resolvent of Brownian motion killed at 0, with $\psi_{\lambda}(x)=\mathbb{E}^{(x,a)}[e^{-\lambda H_{0}}]=e^{- \gamma |x|}$, where $\gamma=\sqrt{2\lambda}$, and $f_{x,a}(y)=f(y,a+y-x)$, we get
\[
\mathcal{V}_{\lambda}f(x,a)=\mathcal{R}_{\lambda}^{-}f_{x,a}(x)+\psi_{\lambda}(x)\mathcal{V}_{\lambda}f(0,a-x)
\]
which follows from the strong Markov property at the hitting time to 0; see e.g. \cite[pg.24]{Howitt} for more details in the one-sided case. Denoting the density of $\mathcal{R}_{\lambda}^{-}$ with respect to Lebesgue measure by $ {r}_{\lambda}^{-}$ yields the desired result.
\end{proof}

Next, we use the joint resolvent of the process $(X_t,B_t)$ derived in Theorem \ref{theo:decompositionXandB} to obtain the conditional distribution at an independent exponential time.
\begin{theorem} \label{thm2}
Let $(X_t,B_t)$  satisfy \eqref{eq:sticky1}-\eqref{eq: sticky2} with initial condition $(0,0)$ and let $T$ be exponentially distributed with rate $\lambda$ independent from $(X,B)$. Then, we have
\[
\mathbb{P}\{X_{T}<x|B_T=b\}= \Big( \frac{\mbox{sgn}(b-x)}{4}+\mbox{sgn}(x) \frac{\gamma}{4\delta}\Big)e^{-\delta|b-x|+\gamma(|b|-{sgn}(x)x)}+\frac{i_x(b)}{2}
\]
where $\gamma^2=2\lambda$, $\delta^2=2\lambda+2\theta\gamma$ and $i_x(b)=1_{\{b<x\}}$.
\end{theorem}
\begin{proof}
 For a bounded Borel function $g$ on $\mathbb{R}^2$, the resolvent $\mathcal{V}_{\lambda}$ of $(X,B)$ satisfies $\mathcal{V}_\lambda g(0,0) =\mathbb{E}^{(0,0)}\,g(X_T,B_T)$. In particular, consider
 \[
 g(y,b)=1_{\{y\le x\}}f(b)=: i_x(y) f(b)
  \]
 for a bounded Borel function $f$ on $\mathbb{R}$. On the other hand, the projection property of conditional expectation implies
\[
\mathbb{E}[g(X_T,B_T)]=\mathbb{E}[1_{\{X_{T}\le x\}}f(B_{T})]=\mathbb{E}[\mathbb{P}\{X_{T} \le x|B_T\}f(B_{T})]\; .
\]
Therefore, if we show that
\[
\mathcal{V}_{\lambda}g(0,0)=[\mathcal{V}_{\lambda}(i_x f)](0,0)= [\mathcal{R}_{\lambda}(e_xf)](0)
 \]
for some $e_x:\mathbb{R}\rightarrow \mathbb{R}$, where $R_\lambda$ is the resolvent of Brownian motion, we can conclude that $e_x(b)=\mathbb{P}\{X_{T}<x|B_T=b\}$ with $(X_0,B_0)=(0,0)$. Now, by Theorem \ref{theo:decompositionXandB}, we get
\begin{eqnarray*}
[\mathcal{V}_\lambda (i_x f)](0,0)&=&\frac{1}{\theta}\int_{-\infty}^{\infty}\nu_\lambda (0,0,0,b)i_x(0)f(b)\, db
\\&&+\int_{-\infty}^{\infty}\int_{-\infty}^{\infty}\nu_\lambda(0,0,y,b)\, i_x(y)f(b)dy\, db\\
&=&
i_x(0)\frac{1}{2\delta}\int_{-\infty}^{\infty}e^{-\delta|b|}f(b)\, db\\
&& +\frac{\theta}{2\delta}\int_{-\infty}^{\infty}\int_{-\infty}^{x}e^{-\gamma|y|
-\delta|b-y|}f(b)\, dy\, db
\end{eqnarray*}
where  $\nu_{\lambda}(x,a,y,b)=\frac{\theta}{2\delta}e^{-\gamma(|y|+|x|)-\delta |b-y-a+x|}$. Assume $x>0$, then
\begin{eqnarray*}
	[\mathcal{V}_\lambda (i_xf)](0,0)&=&\frac{1}{2\delta}\int_{-\infty}^{\infty}e^{-\delta|b|}f(b)db+\frac{\theta}{2\delta}\int_{-\infty}^{\infty}f(b)\int_{-\infty}^{x}e^{-\gamma|y|-\delta|b-y|}dydb\\
	 &=&\frac{1}{2\delta}\int_{-\infty}^{\infty}e^{-\delta|b|}f(b)db+\frac{\theta}{2\delta}\int_{-\infty}^{\infty}f(b)\int_{-\infty}^{x}e^{-\gamma|y|-\delta|b-y|}dydb\\
	&=&\frac{1}{2\delta}\int_{-\infty}^{\infty}e^{-\delta|b|}f(b)db\\
	&&+\frac{\theta}{2\delta}\int_{x}^{\infty}f(b)\Big(\int_{0}^{x}e^{-\gamma y-\delta(b-y)}dy+\int_{-\infty}^{0}e^{\gamma y-\delta(b-y)}dy\Big)db\\
	&& +\frac{\theta}{2\delta}\int_{0}^{x}f(b)\Big(\int_{b}^{x}e^{-\gamma y-\delta(y-b)}dy+\int_{0}^{b}e^{-\gamma y-\delta(b-y)}dy\\
	&&+\int_{-\infty}^{b}e^{\gamma y-\delta(b-y)}dy\Big)db\\
	&&+\frac{\theta}{2\delta}\int_{-\infty}^{0}f(b)\Big(\int_{0}^{x}e^{-\gamma y-\delta(y-b)}dy+\int_{b}^{0}e^{\gamma y-\delta(y-b)}dy\\
	&&+\int_{-\infty}^{b}e^{\gamma y-\delta(b-y)}dy\Big)db\\
	&=:&\int_{-\infty}^{\infty}\frac{1}{\gamma}e^{-\gamma|b|}e^{(1)}_x(b)f(b)db
	\;=\; [\mathcal{R}_{\lambda}(e_x^{(1)}f)](0)
\end{eqnarray*}
where $\mathcal{R}_{\lambda}$ is the resolvent of Brownian motion  and  explicitly
\begin{equation}
e_x^{(1)}(b)=e^{-\delta|b-x|+\gamma(|b|-x)}\Big(\frac{\mbox{sgn}(b-x)}{4}+\frac{\gamma}{4\delta}\Big)+\frac{i_x(b)}{2}.
\end{equation}
For the case $x\leq 0$, we have
\begin{eqnarray*}
[V_{\lambda}(i_{x}f)](0,0)
&=&	 \frac{\theta}{2\delta}\int_{-\infty}^{\infty}f(b)\int_{-\infty}^{x}e^{-\gamma|y|
-\delta|b-y|}dy\,db\\
&=&\frac{\theta}{2\delta}\int_{x}^{\infty}f(b)\int_{-\infty}^{x}e^{\gamma y-\delta(b-y)}dy\,db \\
&&+\frac{\theta}{2\delta}\int_{-\infty}^{x}f(b)\Big(\int_{b}^{x}e^{\gamma y-\delta(y-b)}dy+\int_{-\infty}^{b}e^{\gamma y-\delta(b-y)}dy\Big)\,db\\
&=&\int_{-\infty}^{\infty}\frac{1}{\gamma}e^{-\gamma|b|}e^{(2)}_{x}(b)f(b)db\\
&=&[\mathcal{R}_{\lambda}(e^{(2)}_x f)](0))
\end{eqnarray*}
where
\[ e^{(2)}_x(b)=e^{-\delta|b-x|+\gamma(|b|+x)}\Big(\frac{\mbox{sgn}(b-x)}{4}-\frac{\gamma}{4\delta} \Big)+\frac{i_x(b)}{2}.
\]
Letting $e_x(b)= 1_{\{x>0\}}e^{(1)}_x(b)+1_{\{x\leq 0\}}e^{(2)}_x(b)$,  we  obtain
\[
 \mathbb{E} [1_{\{X_{T}<x\}}f(B_T)]=\mathbb{E} [e_x(B_{T})f(B_T)]	
\]
which completes the proof.
\end{proof}

Theorem \ref{thm2} allows us to calculate the conditional distribution of $X_t$ given $B_t$, for  $t>0$, next.
\begin{corollary}\label{theo: conditionalProbability}
The conditional distribution of sticky Brownian motion $X_t$, given the driving Brownian motion $B_t$, both staring from zero, is given by
\begin{multline}
\mathbb{P}\{X_t\le x | B_t=b\}= \sqrt{2\pi t} \:e^{\frac{b^2}{2t}}\, \Big[\frac{\mbox{sgn}(b-x)}{4}f_1(x,b,t)*f_2(x,b,t)\\
 +\frac{sgn(x)}{4}f_3(x,b,t)*f_4(x,b,t) + \frac{i_x(b)}{2}f_2(x,b,t)\Big] \nonumber
\end{multline}
where $*$ is the convolution operator, and
\begin{eqnarray*}
f_1(x,b,t)&=&  \frac{|b-x|}{\sqrt{2\pi t^3}}e^{-\frac{2|b-x|^2}{4 t}-\theta|b-x|}\\
&&+\frac{\theta}{2\sqrt{2\pi t^3}}\int_{\sqrt{2}|b-x|}^{\infty} \xi \, e^{-\frac{\xi^2}{4t}-\frac{\theta \xi}{\sqrt{2}}}\, I_1\left(\theta\sqrt{\frac{\xi}{2}-|b-x|^2} \right)\, d\xi \\
f_2(x,b,t)&=& \frac{1}{\sqrt{\pi t}}e^{-\frac{|x|^2}{2t}}\\
f_3(x,b,t)&=&\frac{1}{2\sqrt{2\pi t^3}}\int_0^{\infty}\xi e^{ -\frac{\xi^2}{4t}-\frac{\theta \xi}{\sqrt{2}}}\,d\xi \\
f_4(x,b,t)&=&\frac{|x|}{\sqrt{2 \pi t^3}}\, e^{-\frac{\sqrt{2}|x|^2}{4t}}
\end{eqnarray*}
\end{corollary}
\begin{proof}
Under the initial condition $(X_0,B_0)=(0,0)$, we can write the joint distribution as
\begin{eqnarray}
\lefteqn{\mathbb{P}\{X_T\le x, B_T \in db\}=\mathbb{P}\{X_T<x|B_T=b\}\mathbb{P}\{B_T \in db\} } \nonumber\\
& = \displaystyle{\left[\frac{\gamma}{2}\Big( \frac{\mbox{sgn}(b-x)}{4}+\mbox{sgn}(x) \frac{\gamma}{4\delta}\Big)e^{-\delta|b-x|-\gamma |x|  }+\frac{i_x(b)}
 {2}\frac{\gamma}{2}e^{-\gamma|b|}\right]db}  \label{joint}
\end{eqnarray}
where we have used the fact that $B_T$ is double exponentially distributed with parameter $\gamma=\sqrt{2\lambda}$ as $T$ is exponential with parameter $\lambda$. By conditioning on $T$, we have
\begin{eqnarray*}
\mathbb{P}\{X_{T}<x,B_{T} \in db\}&=&\int_0^{\infty}\mathbb{P}\{X_{T}<x,B_{T} \in db | T=t\}\lambda e^{-\lambda t}dt \\
&=& \int_0^{\infty}\mathbb{P}\{X_{t}<x,B_{t} \in db \}\lambda e^{-\lambda t}dt
\end{eqnarray*}
Thus, we can obtain $\mathbb{P}\{X_t<x , B_t \in db \}$ by Laplace inversion of \eqref{joint}.
Denoting the Laplace transform operator by $\mathcal{L}_{\lambda}$, we get
\begin{align*}
\mathbb{P}\{ X_t<x,B_t \in db \}&=\mathcal{L}^{-1}_{\lambda}\Big(\frac{\mbox{sgn}(b-x)}{4\gamma}\exp\{-\delta|b-x|-\gamma |x|\}\Big)\, db\\
&+\mathcal{L}_{\lambda}^{-1}\Big( \frac{1}{4\delta}\mbox{sgn}(x)\exp\{-\delta|b-x|-\gamma |x|\}\Big)\, db\\
&+\mathcal{L}_{\lambda}^{-1}\Big( \frac{i_x(b)}{2\gamma}e^{-\gamma(|b|)}\Big)\, db\\
&= \frac{\mbox{sgn}(b-x)}{4} \mathcal{L}^{-1}_{\lambda}\Big(\exp\{-\delta|b-x|\}\frac{1}{\gamma}\exp\{-\gamma |x|\}\Big)\, db\\
&+\frac{\mbox{sgn}(x)}{4}\mathcal{L}_{\lambda}^{-1}\Big( \frac{1}{\delta}\exp\{-\delta|b-x|\}\exp\{-\gamma |x|\}\Big)\, db\\
&+ \frac{i_x(b)}{2} \mathcal{L}_{\lambda}^{-1}\Big(\frac{1}{\gamma}e^{-\gamma(|b|)}\Big)\, db
\end{align*}
Let the inverse Laplace transforms  $f_i :\mathbb{R}^2 \times [0,\infty)  \rightarrow  \mathbb{R}$ be assigned as
\begin{align*}
&f_1(x,b,t) = \mathcal{L}^{-1}_{\lambda}\Big(\exp\{-\delta|b-x|\}\Big)[t]\, , \quad
f_2(x,b,t) = \mathcal{L}_{\lambda}^{-1}\Big(\frac{1}{\gamma}\exp\{-\gamma |x|\} \Big)[t]\\
&f_3(x,b,t) = \mathcal{L}_{\lambda}^{-1}\Big(\frac{1}{\delta}\exp\{-\delta|b-x|\}\Big)[t]\, , \quad
f_4(x,b,t) = \mathcal{L}_{\lambda}^{-1}\Big(\exp\{-\gamma |x|\} \Big)[t]
\end{align*}
Then,  the inverse Laplace formula in \cite[pp. 978-979]{Handbook}
yields the result.
\end{proof}

\section{Distribution of Occupation Times}

In this section, we find the distributions of occupation times of positive and negative axes by the two-sided sticky Brownian motion. We also find the distribution of the time spent at 0, which we show to be the same as that of the one-sided case. The desired distributions will be calculated via the relation between the two-sided and the one-sided sticky Brownian motions.

We first recover the distribution of the time spent at 0 by the one-sided sticky Brownian motion, which satisfies Equation \eqref{eq:one-sided}. This result is cited in \cite{Warren}, but not given explicitly. For the sake of completeness, we derive the distribution of both occupation times in the following lemma.
\begin{lemma}\label{theo: onesided}
Let $\bar{X}_t$ denote the one-sided sticky Brownian motion satisfying \eqref{eq:one-sided}, and let
\begin{eqnarray*}
\bar{A}_t^0 = \int_{0}^t 1_{\{\bar{X}_s=0\}}ds \quad \mbox{and} \quad
\bar{A}_t^+= \int_{0}^t 1_{\{\bar{X}_s>0\}}ds \; .
\end{eqnarray*}
Then, we have
\begin{eqnarray*}
 i) \; \mathbb{P} \{ \bar{A}^0_s <t \} &=& 1-\text{Erfc}(\frac{\theta t}{\sqrt{2(s-t)}})\\
ii)\; \mathbb{P}\{\bar{A}^+_s <t \}&=& \text{Erfc}(\frac{\theta(s-t)}{\sqrt{2t}})
\end{eqnarray*}
for $t\le s$, where $\text{Erfc}(z):= \frac{2}{\sqrt{\pi }} \int_{z}^{\infty} e^{-\xi^2} d\xi$.
\end{lemma}
\begin{proof}
Let  $\bar{\alpha}_t^+ =\inf \{ s: \bar{A}^+_s \geq t \}$. It has been shown in  \cite{Warren} that $\bar{A}^0_{\bar{\alpha}^+_T}$ is exponentially distributed with rate  $\theta\gamma=\theta \sqrt{2\lambda}$. Formally, we can write
\begin{eqnarray*}
\theta \sqrt{2\lambda} \, e^{-\theta \sqrt{2\lambda}\, x} \,dx =\mathbb{P}\{\bar{A}^0_{\alpha^+_T}\in dx\} &=&\int_0^{\infty}\mathbb{P}\{\bar{A}^0_{\alpha^+_t}\in dx | T=t\}\mathbb{P}\{T\in dt\}\\
&=&\int_0^{\infty} \mathbb{P}\{\bar{A}^0_{\alpha^+_t}\in dx\}\lambda e^{-\lambda t}dt . \;
\end{eqnarray*}
That is,  $(1/\lambda) \mathbb{P}\{\bar{A}^0_{\alpha^+_T}\in dx\}=\mathcal{L}_{\lambda}(f(\cdot,x))$ where $f(t,x) \, dx=\mathbb{P}\{\bar{A}^0_{\alpha^+_t}\in dx\}$. The formula \cite[(10) p.964]{Handbook} yields
\begin{equation}   \label{Warr}
\mathbb{P}\{ \bar{A}^0_{\bar{\alpha}_t^+} \in dx\}=\frac{\theta \sqrt{2}}{\sqrt{\pi t}}e^{-\frac{\theta^2 x^2}{2t}}dx
\end{equation}
for $x\ge 0$, which is half normal distribution with parameter $\sigma = \sqrt{t}/ \theta$.  Recall that $t=\bar{A}^0_t + \bar{A}^+_t$ for each $t>0$.  So,  the distribution of $\bar{\alpha}^{+}_t$ can be calculated  as
\begin{eqnarray*}
\mathbb{P}\{\bar{\alpha}_t^+<x\}&=&\mathbb{P}\{t+\bar{A}^0_{\bar{\alpha}_t^+}<x\} \nonumber \; =\;\mathbb{P}\{\bar{A}^0_{\bar{\alpha}_t^+}<x-t\} \\
& = &\int_0^{x-t}\frac{\theta \sqrt{2}}{\sqrt{\pi t}}e^{-\frac{\theta^2 y^2}{2t}}dy \label{eq1}
\;=\; \int_t^{x}\frac{\theta \sqrt{2}}{\sqrt{\pi t}}e^{-\frac{\theta^2 (y-t)^2}{2t}}dy
\end{eqnarray*}
in view of \eqref{Warr}. Then, the density is given by
\[
\mathbb{P}\{\bar{\alpha}_t^+ \in dy\}=\frac{\theta \sqrt{2}}{\sqrt{\pi t}}e^{-\frac{\theta^2 (y-t)^2}{2t}}dy  \quad \quad \mbox{for} \quad y \ge t
\]
which is also half Normal, but shifted to $t$.  The result i) follows from the definition of $\bar{\alpha}_t^{+}$ since
\begin{equation}
    \mathbb{P}\{\bar{A}^+_s<t\}=\mathbb{P}\{s<\bar{\alpha}^{+}_t\}
\end{equation}
Then, ii) follows from the fact that $\mathbb{P}\{
\bar{A}^0_s<t\}=\mathbb{P}\{s-\bar{A}^+_s<t\}=1-\mathbb{P}\{\bar{A}^+_s\leq (s-t)\}$.
\end{proof}

Recall that the occupation times of the sticky Brownian motion have been defined in equations \eqref{eq: zerotime} to \eqref{eq:negativetime}. Clearly, $A^+$ and $A^-$ have the same distribution,  which is given in the following  main result of this section. The distribution   deviates from the Arcsine law corresponding to Brownian motion \cite[Thm.VI.2.7]{RevuzYor}. Although our method of proof is borrowed from that for Brownian motion, we take into account that $A_t^0$ is positive with positive probability in the sticky case.

\begin{theorem} \label{thm3} Let $X$ be the sticky Brownian motion satisfying \eqref{eq:sticky1}-\eqref{eq: sticky2} with $X_0=0$. For $0<t<1$, the occupation times of $X$ satisfy
\begin{eqnarray*}
i) \; \mathbb{P}\{A^0_1>t\}&=&\text{Erfc}(\frac{\theta t}{\sqrt{2(1-t)}})  \\
ii)\;  \mathbb{P}\{A^+_1 > t \}&=&  1-\mathbb{E}\, \left[ \text{Erfc}
\left( \frac{ -Z^2+Z\sqrt{Z^2+\theta^2(1-t)}} {\theta \sqrt{2t}}\right) \right]
\end{eqnarray*}
for $0<t<1$, where $Z$ is a standard Gaussian random variable.
\end{theorem}

\begin{proof} i)
The relation between the two-sided and one-sided sticky Brownian motion has   been given in \cite{Engelbert}. Through the specific construction of $(X_t, B_t)$ given in \cite[Sec.2]{Engelbert} satisfying \eqref{eq:sticky1}-\eqref{eq: sticky2}, the pair $(|X_t|, \hat{B}_t)$ satisfies Equation \eqref{eq:one-sided} for the one-sided sticky Brownian motion and equivalently the system
\begin{align}
dX_{t}&=1_{\{X_{t}\neq 0\}} dB_t + \frac{1}{2}1_{\{X_t>0\}}dB_t\label{eq:oneside1}\\ 
1_{\{X_t=0\}}dt&= \frac{1}{2\theta}dl_t^0(X)\label{eq:oneside2}
\end{align}
by \cite[Thm.5]{Engelbert}.
Note that local time of two-sided sticky Brownian motion $X_t$ satisfies
\begin{align*} 
l_t^0(X)&=\mathbb{P}\textit{-}\lim_{\varepsilon \downarrow0}\frac{1}{2\varepsilon}\int_0^tI(-\varepsilon \leq X_s\leq \varepsilon ) d\langle X,X\rangle_s \\
&=\mathbb{P}\textit{-}\lim_{\varepsilon \downarrow0}\frac{1}{2\varepsilon}\int_0^tI(0 \leq |X_s|\leq \varepsilon) d\langle X,X\rangle_s, 
\end{align*}
and local time of one-sided sticky Brownian motion $|X_t|$ is defined as $l_0^t(|X|)=\mathbb{P}\textit{-}\lim_{\varepsilon \downarrow 0}\frac{1}{\varepsilon}\int_0^t I(0\leq |X_s| \leq \varepsilon)d\langle X_s,X_s\rangle$ which implies $l_t^0(X)=\frac{1}{2}l_t^0(|X|)$. In view of SDE's \eqref{eq:one-sided} (equivalently \eqref{eq:oneside1}-\eqref{eq:oneside2}) and \eqref{eq:sticky1}-\eqref{eq: sticky2}, we obtain 
\[
A^0_t=\frac{1}{\theta} l_t^0(X)=\frac{1}{2 \theta}l_t^0(|X|)=\bar{A}_t^0,
\]
as in \cite[Eq.3.8]{Engelbert}, where $\bar{A}_t^0$ is as defined in Lemma \ref{theo: onesided}. This implies that $A^0_t$ and $\bar{A}^0_{t}$ have the same distribution that is also given in Lemma \ref{theo: onesided}.

ii) Note that $t=A^0_t+A^-_t+A^+_t$, which implies
\begin{equation} \label{identity}
\alpha^+_t=t+ A^0_{\alpha^+_t}+A^{-}_{\alpha^+_t}\; .
 \end{equation}
Since $\mathbb{P}\{A^+_1>t\}=\mathbb{P}\{1>\alpha^+_t\}$, we will  study $\alpha^+_t$ to get the distribution of $A^+_1$.  Recall from Lemma \ref{W} that $\frac{\theta}{2}A^0_{\alpha^+_t}$ is the running supremum of the Brownian motion $-W^+$. Hence, we have
\begin{equation} \label{first}
A^{0}_{\alpha_{t}^{+}}=\frac{2}{\theta}\: \underset{s\leq t}{\sup}(-W^{+}_{s})=:\frac{2}{\theta}\, S^+_t
\end{equation}
at time $t$, where $S^+_t $ has the law of a reflected Brownian motion. Now, for studying the  term $A^-_{\alpha^+_t}$ in \eqref{identity}, define the stopping times
\[
\tau_t=\inf\{s>0: l^0_s(X)>t\}
 \quad \mbox{and} \quad  T^{W^-}(a)=\inf \{t : W^-_t >a \}
\]
for $a>0$. Note that \eqref{eq:Tanaka2} implies
\[
0=X^-_{\tau_t}=-W^-_{A^-_{\tau_t}}+\frac{\theta}{2}A^0_{\tau_t}
=-W^-_{A^-_{\tau_t}}+\frac{l^0_{\tau_t}(X)}{2}
=-W^{-}_{A^-_{\tau_t}}+\frac{t}{2}
\]
in view of \eqref{eq: sticky2}. This yields $W^-_{A^-_{\tau_t}}=\frac{t }{2}$. Then, we have
\[
A^-_{\tau_t}=T^{W^-}(t/2)
\]
for each $t\ge 0$, which follows along the same lines of proof for Brownian motion in \cite[Thm.VI.2.7]{RevuzYor}. Replacing $t$ with $l^0_{\alpha^+_t}$ above yields
\begin{equation}   \label{arcsine}
A^-_{\alpha^+_{t}} {=} T^{W^{-}}(l^0_{\alpha^+_t}/2) = T^{W^{-}}\left (\frac{\theta}{2} A^0_{\alpha^+_t} \right) = T^{W^{-}} ({S^+_t})
\end{equation}
where the last equality follows from \eqref{first}.

 We are now ready to compute $\mathbb{P}\{A^+_1 > t \}$. Using \eqref{identity}, \eqref{first} and \eqref{arcsine}, we get
\begin{eqnarray}
\mathbb{P}\{A^+_1 > t \}&=&\mathbb{P} \{1> \alpha^+_t\}= \mathbb{P}\{0> t-1+A^0_{\alpha^+_t}+A^-_{\alpha^+_t}\} \nonumber \\
&=& \mathbb{P} \{ 0>t-1+ \frac{2}{\theta}S^+_t+ T^{W^{-}} ({S^+_t}) \} \nonumber \\
&=& \mathbb{E} [\mathbb{P} \{ 0>t-1+ \frac{2}{\theta}S^+_t+ T^{W^{-}} ({S^+_t})\,|\, S^+_t  \}  ] \nonumber \\
&=& \mathbb{E} [\mathbb{P} \{ 0>t-1+ \frac{2}{\theta}S^+_t+ \frac{(S^+_t)^2}{(W^-_1)^2}\,|\,S^+_t  \}  ] \label{theta}\\
&=& \mathbb{P} \left\{ 0>t-1+ \frac{2}{\theta}S^+_t+ \frac{(S^+_t)^2}{(W^-_1)^2} \right \}\nonumber
\end{eqnarray}
where \eqref{theta} follows by the independence of $S^+$ and $W^-$ from \cite[III.Prop.3.10]{RevuzYor}, which states that
\[
T^{W^{-}} (a)  \stackrel{(d)}{=} \frac{a^2}{(W_1^-)^2} \quad \quad a>0 \; .
\]
As   $S^+$ and $W^-$ are independent, we  will only need the roots of the polynomial $P(x)=t-1+\frac{2}{\theta}x+\frac{x^2}{w^2}$. The positive root is given by $(w^2/\theta)(-1+\sqrt{1+\theta^2(1-t)/w^2}) $. As a result, we get
\begin{eqnarray*}
\mathbb{P}\{A^+_1 > t \}&= &\int_{-\infty}^{\infty} \left(\int_{0}^{(x^2/\theta)(-1+\sqrt{1+\theta^2(1-t)/x^2})}  \frac{2 \,e^{-y^2/(2t)}}{\sqrt{2\pi t}}  \, dy \right)\: \frac{e^{-x^2/2}}{\sqrt{2\pi}}\, dx \\
&=& \int_{-\infty}^{\infty}  \left( \int_{0}^{\frac{x^2}{\theta \sqrt{t}} (-1+\sqrt{1+\theta^2(1-t)/x^2})}  \frac{2\, e^{-y^2/2}}{\sqrt{2\pi}} \, dy \right)  \frac{e^{-x^2/2}}{\sqrt{2\pi}} \,  dx
\end{eqnarray*}
since $S^+$ has the law of a reflected Brownian motion.
\end{proof}

As a final remark, we observe  from \eqref{theta} that the distribution of $A_1^+$ converges to the Arcsine law as $\theta \rightarrow \infty$, by similar computations that would follow as in the proof of \cite[Thm.VI.2.7]{RevuzYor} for Brownian motion. This is consistent with the fact that the two-sided sticky Brownian motion approximates a standard Brownian motion in this case.

\vspace{5mm}
\noindent \textbf{Acknowledgement.} This research is supported by TUBITAK Project 115F086. The authors would like to thank the anonymous reviewer for valuable comments that  improved the manuscript.

\bibliographystyle{plain}


\end{document}